\documentclass[12pt, twoside, reqno]{article}
\usepackage{amsmath,amsthm}
\usepackage{amssymb}
\usepackage{enumitem}

\newtheorem{theorem}{Theorem}[section]

\newtheorem{lemma}[theorem]{Lemma}

\newtheorem{problem}[theorem]{Problem}

\theoremstyle{definition}

\numberwithin{equation}{section}

\begin{document}

\baselineskip=17pt


\title{On the sum of a prime and two Fibonacci numbers}

\author{Ji-Zhen Xu$^{a,b}$ and  Yong-Gao Chen$^{a}$\footnote{Corresponding author} \\
\small $^a$School of Mathematical Sciences and Institute of Mathematics,\\
\small Nanjing Normal University,  Nanjing  210023,  P.R. China\\
\small $^b$Nanjing Vocational College of Information Technology,\\
 \small Nanjing 210023, P.R. China\\
\small E-mails:965165607@qq.com (Xu), ygchen@njnu.edu.cn (Chen)  }
\date{}

\maketitle


\begin{abstract}
 Let $\{f_n\}$ be
the  Fibonacci sequence. For any positive integer $n$, let $r(n)$
be the number of solutions of $n=p+f_{k_1^{2}} +f_{k_{2}^{2}}$,
where $p$ is a prime and $k_1, k_2$
 are nonnegative integers with $k_1\le k_2$.
In this paper, it is proved that $\{ n : r(n)=0\} $ contains an
infinite arithmetic progression, and both sets $\{ n : r(n)=1\} $
and $\{ n : r(n)\ge 2\}$ have  positive asymptotic densities.
\end{abstract}

\renewcommand{\thefootnote}{}



\renewcommand{\thefootnote}{\arabic{footnote}}
\setcounter{footnote}{0}

\section{Introduction}

In 1934, Romanoff \cite{Romanoff1934} proved that there is a
positive proportion of positive odd numbers that can be
represented as the sum of a prime and a power of 2. In 1950, Erd\H
os \cite{Erdos1950} proved that there is an infinite arithmetic
progression of positive odd numbers none of which can be
represented as  the sum of a prime and a power of 2. For related
progress, one may refer to Crocker \cite{Crocker1971},  Elsholtz
and Schlage-Puchta \cite{Elsholtz2018}, Pan \cite{Pan2011}, Pintz
\cite{Pintz2006} and Xu and Chen \cite{XuChenIJNT2024}.

Let $\mathcal{P}$ be the set of all positive primes, $\mathbb{N}$
the set of all positive integers and $\mathbb{N}_0=\mathbb{N}\cup
\{ 0\}$.  Let $\{f_n\}$ be the  Fibonacci sequence which is
defined by $f_0=0$, $f_1=1$ and $f_n=f_{n-1}+f_{n-2}$ $(n\ge 2)$.
In 2010, Lee \cite{Lee} proved that the set of positive integers
that can be represented as $p+f_{k}$, $p\in \mathcal{P}$, $k\in
\mathbb{N}$, has positive lower asymptotic density. In 2021, Liu
and Xue \cite{Liu} showed that this lower asymptotic density is
more than $0.025$.  Wang and Chen \cite{WangChen2023} improved
$0.025$ to $0.143$. Chen and Wang \cite{ChenWang2024}  showed that
there is a positive proportion of positive integers that can be
represented uniquely as $p+f_{k}$, $p\in \mathcal{P}$, $k\in
\mathbb{N}$.
 Recently, Chen and
Xu \cite{ChenXu2024} proved that if $r_1, \dots, r_t$ are positive
integers such that $r_1^{-1}+ \dots +r_t^{-1}\ge 1$, then  there
is a positive proportion of positive integers which can be
represented as $p+f_{k_1^{r_1}}+\cdots +f_{k_{t}^{r_t}}$, where
$p\in \mathcal{P}$ and $k_1,\cdots ,  k_t\in \mathbb{N}$. On the
other hand,  $\check{S}$iurys \cite{Siurys} showed that there is
an infinite arithmetic progression of positive integers none of
which can be represented as $p+f_{k}$, where $p\in \mathcal{P}$
and $k\in \mathbb{N}_0$.  Chen and Xu \cite{ChenXu2024}
conjectured that if $r_1, \dots, r_t$ are positive integers such
that $r_1^{-1}+ \dots +r_t^{-1}= 1$, then there is an infinite
arithmetic progression of positive integers none of which cannot
be represented as $p+f_{k_1^{r_1}}+\cdots +f_{k_{t}^{r_t}}$, where
$p\in \mathcal{P}$ and $k_1,\cdots , k_t\in \mathbb{N}$.

In this paper, we concern the integers of the form  $p+f_{k_1^{2}}
+f_{k_{2}^{2}}$, where $p\in \mathcal{P}$ and $k_1, k_2\in
\mathbb{N}_0$. The following results are proved.

\begin{theorem}\label{thm1}
There is an infinite arithmetic progression of positive integers
none of which can be represented as
 $p+f_{k_1^{2}} +f_{k_{2}^{2}}$, where $p\in \mathcal{P}$ and $k_1, k_2\in
\mathbb{N}_0$.
\end{theorem}

\begin{theorem}\label{thm3}
There is a positive proportion of positive integers which can be
uniquely represented as $p+f_{k_1^{2}}+f_{k_2^{2}}$, where $p\in
\mathcal{P}$ and $k_1, k_2\in \mathbb{N}_0$ with $k_1\leq k_2$.
\end{theorem}

\begin{theorem}\label{thm4}
There is a positive proportion of positive integers which can be
represented in at least two ways as $p+f_{k_1^{2}}+f_{k_2^{2}}$,
where $p\in \mathcal{P}$ and $k_1, k_2\in \mathbb{N}_0$ with
$k_1\leq k_2$.
\end{theorem}

Similar to the proof of Theorem \ref{thm4}, one can prove that
there is a positive proportion of positive integers which can be
represented in at least three ways as $p+f_{k}$, where $p\in
\mathcal{P}$ and $k\in \mathbb{N}_0$. Let
$$r(n)=\sharp \{ (p, k_1, k_2 ): n=p+f_{k_1^2}+f_{k_2^2}, p\in \mathcal{P}, k_1, k_2\in \mathbb{N}_0, k_1\le
k_2\} .$$ By Chen and Ding \cite{ChenDing2022} or
\cite{ChenDing2023}, we have
$$\limsup_{n\to \infty} \frac{r(n)}{\log\log n} >0.$$
The initial positive integers, which cannot be represented as
 $p+f_{k_1^{2}} +f_{k_{2}^{2}}$, where $p\in \mathcal{P}$ and $k_1, k_2\in
\mathbb{N}_0$, are $1,28, 122, 125, 156, 178, 189, 190, 206, \dots
$. We pose the following three problems for further research.

\begin{problem}\label{prob1} Is there a positive integer $r$
such that every integer greater than $1$ can be written as the
form $p+f_{k_1^{2}}+f_{k_2^{2}}+\cdots +f_{k_r^{2}}$? where $p\in
\mathcal{P}$ and $k_i\in \mathbb{N}_0$ $(1\le i\le r)$.
\end{problem}

\begin{problem}\label{prob2} Is it true that for any positive integer $m$, there is a positive proportion of positive integers which can be
represented in at least $m$ ways as $p+f_{k_1^{2}}+f_{k_2^{2}}$?
where $p\in \mathcal{P}$ and $k_1, k_2\in \mathbb{N}_0$ with
$k_1\leq k_2$.
\end{problem}

\begin{problem}\label{prob3} Is it true that for any positive integer $m$, there is a positive proportion of positive integers which can be
represented in at least $m$ ways as $p+f_k$? where $p\in
\mathcal{P}$ and $k\in \mathbb{N}_0$.
\end{problem}

 Problem \ref{prob1} is verified for $r=6$ and $1<n\le 10^8$.
 For $r=5$, the least three counterexamples are $12877723$,
 $13445485$ and $14542811$. So $r\ge 6$. Generally, we have a similar problem for
 $p+f_{k_1^s}+\cdots +f_{k_r^s}$, where $s$ is a positive integer.
 Elsholtz and Schlage-Puchta
\cite{Elsholtz2018} proved that for any positive integer $m$,
there is a positive proportion of positive integers which can be
represented in at least $m$ ways as $p+2^k$, where $p\in
\mathcal{P}$ and $k\in \mathbb{N}$. Currently, we cannot prove or
disprove that  for any positive integer $m$, there is a positive
proportion of positive integers which can be represented in at
least $m$ ways as $p+2^{k_1^{2}}+2^{k_2^{2}}$, where $p\in
\mathcal{P}$ and $k_1, k_2\in \mathbb{N}_0$ with $k_1\leq k_2$.

 In this paper, let
$a\hskip -1mm\pmod{m} =\{ a+mk : k=0, \pm 1, \dots \} $. For any
prime $q$, let $u(q)$ be the least positive integer $u$ such that
$f_u\equiv f_0\hskip -1mm\pmod{q}$ and $f_{u+1}\equiv f_1\hskip
-1mm\pmod{q}$. It follows from the definition of the Fibonacci
numbers that $f_{n+u(q)}\equiv f_n\hskip -1mm\pmod{q}$ for all
nonnegative integers $n$. For any positive integer $d>1$, we
define $v(d)=\max\{u(p): p\mid d, p\in\mathcal{P}\}$. We appoint
$v(1)=1$.

\section{Proof of Theorem \ref{thm1} }\label{sec1}

In this section, we prove the following result which will be used
to prove Theorems \ref{thm1} and \ref{thm3}. It is also
interesting itself. Let
\begin{eqnarray*}  M=621386267972593776074029725204132260351094,\end{eqnarray*}
\begin{eqnarray*}n_0=145090702762087107496178817377891453382310, \end{eqnarray*}
\begin{eqnarray*}p_0=612271766655191184669165916436108277668671,\end{eqnarray*}
\begin{eqnarray*}\mathcal{P}_M=\{ p : p \in \mathcal{P}, p\equiv
p_0\hskip -3mm\pmod{M} \} ,\end{eqnarray*}
\begin{eqnarray*}\mathcal{K}_1=\{ k : k\in \mathbb{N}_0, k\equiv 0\hskip -3mm\pmod{192} \} , \end{eqnarray*}
\begin{eqnarray*}\mathcal{K}_2=\{ k :  k\in \mathbb{N}_0, k\equiv 64, 128\hskip -3mm\pmod{192}\}
.\end{eqnarray*} Then
\begin{eqnarray*}M&=&2\times 3\times 7\times 17\times 19\times 47\times 107\times 127\times
 1087\\ &&\times 2207\times 4481\times 6143\times 12289\times 21503\times 34303\times 119809.\end{eqnarray*}

\begin{theorem}\label{thm2}
Let $n$ be a positive integer. Suppose that $n$ can be represented
as $p+f_{k_1^{2}}+f_{k_2^{2}}$, where $p\in \mathcal{P}, k_1, k_2
\in \mathbb{N}_0$. Then $n\equiv n_0\hskip -1mm\pmod{M} $ if and
only if $p\in \mathcal{P}_M$, $k_i\in \mathcal{K}_i$  $(i=1,2)$ or
$k_i\in \mathcal{K}_{3-i}$ $(i=1,2)$.
\end{theorem}

\begin{proof}
By the Chinese remainder theorem, we have
\begin{eqnarray}\label{eq0}
&  & 0\hskip -3mm\pmod{2} \cap 1\hskip -3mm\pmod{3}\cap 3\hskip
-3mm\pmod{7}\cap 1\hskip -3mm\pmod{17}\cap 9\hskip
-3mm\pmod{19}\nonumber
\\
&&\cap~ 6\hskip -3mm\pmod{47}\cap 15\hskip -3mm\pmod{107}\cap 80
\hskip -3mm\pmod{127}\cap
 887\hskip -3mm\pmod{1087}\nonumber
\\
&&\cap~ 987\hskip -3mm\pmod{2207}\cap 58\hskip -3mm\pmod{4481}\cap 3107\hskip -3mm\pmod{6143}\cap 5881\hskip -3mm\pmod{12289} \nonumber
\\
&&\cap~
15722\hskip -3mm\pmod{21503}\cap 2489 \hskip
-3mm\pmod{34303}\cap 115550\hskip -3mm\pmod{119809}\nonumber\\
&=& n_0\hskip -3mm\pmod{M}
\end{eqnarray}
and
\begin{eqnarray}\label{eq1}
&& 1\hskip -3mm\pmod{2} \cap 1\hskip -3mm\pmod{3}\cap 3\hskip
-3mm\pmod{7}
\cap 15 \hskip -3mm\pmod{17} \cap 10\hskip -3mm\pmod{19} \nonumber\\
&& \cap~ 6 \hskip -3mm\pmod{47}\cap 36\hskip -3mm\pmod{107}\cap 80
\hskip -3mm\pmod{127} \cap 887\hskip -3mm\pmod{1087}\nonumber
\\
&& \cap~ 987\hskip -3mm\pmod{2207}\cap 58\hskip
-3mm\pmod{4481}\cap 4227\hskip -3mm\pmod{6143}\cap 11762\hskip -3mm\pmod{12289} \nonumber
\\
&&\cap ~15722\hskip -3mm\pmod{21503}\cap 2489 \hskip
-3mm\pmod{34303}\cap 115550\hskip -3mm\pmod{119809}\nonumber\\
&=&p_0\hskip -3mm\pmod{M}.
\end{eqnarray}

 A
calculation shows that
\begin{eqnarray}\label{e4}
&u(2)=3, u(3)=2^3, u(7)=2^4, u(17)=2^2\times 3^2, u(19)=2\times 3^2,\nonumber\\
& u(47)=2^5,u(107)=2^3\times 3^2, u(127)=2^8, u(1087)=2^7,
\nonumber\\
& u(2207)=2^6, u(4481)=2^7, u(6143)=3\times 2^{12}, u(12289)=3\times 2^{11},\nonumber\\
&  u(21503)=2^{11}, u(34303)=2^{10}, u(119809)=2^9.
\end{eqnarray}

Suppose that $n\equiv n_0\hskip -1mm\pmod{M}$. We will prove that
$k_i\in \mathcal{K}_i$  $(i=1,2)$ or $k_i\in \mathcal{K}_{3-i}$
$(i=1,2)$.

Let $$k_1\equiv a\hskip -3mm\pmod{128},\quad k_2\equiv b\hskip
-3mm\pmod{128}, -64\leq a, b \leq 64. $$ Then $k_1^2\equiv
a^2\hskip -1mm\pmod{128}$ and $k_2^2\equiv b^2\hskip
-1mm\pmod{128}$. Since $u(4481)=128$ and $4481\mid M$, it follows
from  \eqref{eq0} that
\begin{eqnarray}\label{e4a}p=n-f_{k_1^{2}}-f_{k_2^{2}}\equiv n_0-f_{a^2}- f_{b^2}\equiv 58-f_{a^2}- f_{b^2}\hskip -3mm\pmod{4481}.\end{eqnarray}
A calculation shows that \eqref{e4a} implies
\begin{eqnarray}\label{e2}
 p\notin
 &\{2, 3, 7, 17, 19,  47, 107, 127, 1087, 2207, \nonumber\\
 &6143, 12289, 21503, 34303, 119809 \}.
\end{eqnarray}

Let $i\in \{ 1, 2\} $. If $k_i>0$, then there are two nonnegative
integers $\alpha_i$ and $\ell_i$ with $2\nmid \ell_i$ such that
$k_i=2^{\alpha_i} \ell_i $. If $k_i=0$, then let $\alpha_i=+\infty
$ and $\ell_i=0$. Without loss of generality, we may assume that
$\alpha_1\le \alpha_2$.

We divide into the following cases:

{\bf Case 1:} $3\nmid k_1, 3\nmid k_2$ or $3|k_1, 3|k_2$. Then
$k_1^2\equiv k_2^2\equiv 1\hskip -1mm\pmod{3}$ or $k_1^2\equiv
k_2^2 \equiv 0\hskip -1mm\pmod{3}$. Since $u(2)=3$, we have
$f_{k_1^{2}}+f_{k_2^{2}}\equiv 0\hskip -1mm\pmod{2}$. It follows
from  \eqref{eq0} that $p= n-f_{k_1^{2}}-f_{k_2^{2}}\equiv 0\hskip
-1mm\pmod{2}$. So $p=2$, a contradiction with \eqref{e2}. Hence,
either $3\mid k_1$, $3\nmid k_2$ or $3\nmid k_1$, $3\mid k_2$.

{\bf Case 2:} $\alpha_1= \alpha_2=0$. Without loss of generality,
we may assume that $3\mid k_1$ and $3\nmid k_2$. Thus,
$k_1^2\equiv 9\hskip -1mm\pmod{72}$ and $k_2^2\equiv 1, 25,
49\hskip -1mm\pmod{72}$. For $q\in \{17, 19, 107\} $, we have
$u(q)\mid 72$. If $k_2^2\equiv 1\hskip -1mm\pmod{72}$, then by
\eqref{eq0},
$$p=n-f_{k_1^{2}}-f_{k_2^{2}}\equiv 1-f_9-f_1\equiv 0\hskip -3mm\pmod{17}.$$
If $k_2^2\equiv 25\hskip -1mm\pmod{72}$, then by \eqref{eq0},
$$p=n-f_{k_1^{2}}-f_{k_2^{2}}\equiv 9-f_9-f_{25}\equiv 0\hskip -3mm\pmod{19}.$$
If $k_2^2\equiv 49\hskip -1mm\pmod{72}$, then by \eqref{eq0},
$$p=n-f_{k_1^{2}}-f_{k_2^{2}}\equiv 15-f_9-f_{49}\equiv 0\hskip -3mm\pmod{107}.$$
So $p\in \{ 17, 19, 107\}$, a contradiction with \eqref{e2}.

{\bf Case 3:} $\alpha_1=0$ and $\alpha_2>0$. Then $k_1^2\equiv
1\hskip -1mm\pmod{8}$  and $k_2^2\equiv 0,4\hskip -1mm\pmod{8}$.
 Since $u(3)=8$ and $f_4=3$, it follows from  \eqref{eq0} that
 $$p=n-f_{k_1^{2}}-f_{k_2^{2}}\equiv n-1\equiv 0\hskip -3mm\pmod{3}.$$ So
 $p=3$, a contradiction with \eqref{e2}.

{\bf Case 4:}  $1\le \alpha_1\le 4$. Let
$$p_{1,1}=47, \quad p_{1, 1+}=7, \quad p_{2,2}=1087,\quad  p_{2,
2+}=2207,$$ $$ p_{3,3}=119809,\quad  p_{3, 3+}=127, \quad
p_{4,4}=21503,\quad p_{4, 4+}=34303.$$

If $\alpha_1=\alpha_2$, then
 $k_i^2\equiv 2^{2\alpha_1}\hskip -1mm\pmod{2^{2\alpha_1+3}}$ $(i=1,2)$.
Since $u(p_{\alpha_1, \alpha_1})=2^{2\alpha_1+3}$, it follows from
\eqref{eq0} that
 $$p=n-f_{k_1^{2}}-f_{k_2^{2}}\equiv n-2f_{2^{2\alpha_1}}\equiv 0\hskip -3mm\pmod{p_{\alpha_1, \alpha_1}}.$$ So
 $p=p_{\alpha_1, \alpha_1}$, a contradiction with \eqref{e2}.

If $\alpha_1<\alpha_2$, then
 $k_1^2\equiv 2^{2\alpha_1}\hskip -1mm\pmod{2^{2\alpha_1+3}}$ and
 $k_2^2\equiv 0\hskip -1mm\pmod{2^{2\alpha_1+2}}$.
Since $u(p_{\alpha_1, \alpha_1+})=2^{2\alpha_1+2}$, it follows
from  \eqref{eq0} that
 $$p=n-f_{k_1^{2}}-f_{k_2^{2}}\equiv n-f_{2^{2\alpha_1}}\equiv 0\hskip -3mm \pmod{p_{\alpha_1, \alpha_1+}}.$$ So
 $p=p_{\alpha_1, \alpha_1+}$, a contradiction with \eqref{e2}.

 {\bf Case 5:}  $\alpha_1=\alpha_2=5$. Without loss of generality, we may assume that $3\nmid k_1$ and $3\mid k_2$.
 Then $$k_1^2\equiv 2^{10}\hskip -3mm\pmod{3\times 2^{13}}, \quad k_2^2\equiv 2^{10}+2^{13}\hskip -3mm\pmod{3\times 2^{13}}.$$
 Since $u(12289)=3\times 2^{11}$, we have
 $$f_{k_1^2}+f_{k_2^2}\equiv f_{2^{10}}+f_{2^{10}+2^{13}}\equiv 5881\hskip -3mm\pmod{12289}.$$
 It follows that
 $$p=n-f_{k_1^2}-f_{k_2^2}\equiv n-5881\equiv 0\hskip -3mm\pmod{12289}.$$
So $p=12289$, a contradiction with \eqref{e2}.

 {\bf Case 6:}  $\alpha_1=5$ and $\alpha_2>5$. If $3\nmid k_1$ and $3\mid k_2$, then
 $$k_1^2\equiv 2^{10}\hskip -3mm\pmod{3\times 2^{12}}, \quad k_2^2\equiv 0\hskip -3mm\pmod{3\times 2^{12}}.$$
 Since $u(12289)=3\times 2^{11}$, we have
 $$f_{k_1^2}+f_{k_2^2}\equiv f_{2^{10}}+f_0\equiv 5881\hskip -3mm\pmod{12289}.$$
 It follows that
 $$p=n-f_{k_1^2}-f_{k_2^2}\equiv n-5881\equiv 0\hskip -3mm\pmod{12289}.$$
So $p=12289$, a contradiction with \eqref{e2}.

If $3\mid k_1$ and $3\nmid k_2$, then
 $$k_1^2\equiv 2^{10}+2^{13}\hskip -3mm\pmod{3\times 2^{12}}, \quad k_2^2\equiv 2^{12}\hskip -3mm\pmod{3\times 2^{12}}.$$
 Since $u(6143)=3\times 2^{12}$, we have
 $$f_{k_1^2}+f_{k_2^2}\equiv f_{2^{10}+2^{13}}+f_{2^{12}}\equiv 3107\hskip -3mm\pmod{6143}.$$
 It follows that
 $$p=n-f_{k_1^2}-f_{k_2^2}\equiv n-3107\equiv 0\hskip -3mm\pmod{6143}.$$
So $p=6143$, a contradiction with \eqref{e2}.

Up to now, we have proved that if $n\equiv n_0\hskip
-1mm\pmod{M}$, then $\alpha_2\ge\alpha_1\ge 6$ and either $3\mid
k_1$, $3\nmid k_2$ or $3\nmid k_1$, $3\mid k_2$. That is, $k_i\in
\mathcal{K}_i$  $(i=1,2)$ or $k_i\in \mathcal{K}_{3-i}$ $(i=1,2)$.

To prove Theorem \ref{thm2}, we may always assume that $k_i\in
\mathcal{K}_i$  $(i=1,2)$ or $k_i\in \mathcal{K}_{3-i}$ $(i=1,2)$.
Without loss of generality, we may assume that $k_i\in
\mathcal{K}_i$  $(i=1,2)$.  We will prove that $n\equiv n_0\hskip
-1mm\pmod M$ if and only if $p\equiv p_0\hskip -1mm\pmod M$.
Noting that $n=p+f_{k_1^{2}}+f_{k_2^{2}}$, it is enough to prove
that
\begin{equation*}f_{k_1^{2}}+f_{k_2^{2}}\equiv
n_0-p_0\hskip -3mm\pmod M. \end{equation*} That is, for every
prime
\begin{eqnarray*}
 q\in
 &\{2, 3, 7, 17, 19,  47, 107, 127, 1087, 2207, 4481,\\
 &6143, 12289, 21503, 34303, 119809 \} ,
\end{eqnarray*}
 we have
\begin{equation}\label{eq2}f_{k_1^{2}}+f_{k_2^{2}}\equiv
n_0-p_0\hskip -3mm\pmod q.\end{equation}

 We divide into the following cases:

{\bf Case I:} $q=2$. Since $u(2)=3$, $k_1^2\equiv 0\hskip
-1mm\pmod{3}$ and $ k_2^2 \equiv 1\hskip -1mm\pmod{3}$, it follows
from \eqref{eq0} and \eqref{eq1} that \eqref{eq2} holds.

{\bf Case II:} $q\in \{ 17, 19, 107\} $. Then $u(q)\mid 2^3\times
3^2$. By $k_i\in \mathcal{K}_i$ $(i=1,2)$, we have $k_1^2\equiv
0\hskip -1mm\pmod{72}$ and $k_2^2\equiv 16, 40, 64\hskip
-1mm\pmod{72}$.

If $k_2^2\equiv 16\hskip -1mm\pmod{72}$, then
$$f_{k_1^{2}}+f_{k_2^{2}}\equiv f_{16}\equiv 18\equiv n_0-p_0\hskip -3mm\pmod{19}.$$

If $k_2^2\equiv 40\hskip -1mm\pmod{72}$, then
$$f_{k_1^{2}}+f_{k_2^{2}}\equiv f_{40}\equiv 3\equiv
n_0-p_0\hskip -3mm\pmod{17}.$$

If $k_2^2\equiv 64\hskip -1mm\pmod{72}$, then
$$f_{k_1^{2}}+f_{k_2^{2}}\equiv f_{64}\equiv 86\equiv
n_0-p_0\hskip -3mm\pmod{107}.$$

{\bf Case III:} $q\in \{ 3, 7, 47, 127, 1087, 2207, 4481, 21503, 34303,
119809\}$. Then $u(q)\mid 2^{11}$. Thus, $u(q)\mid k_i^2$ by
$2^6\mid k_i$ $(i=1,2)$. It follows that
$$f_{k_1^{2}}+f_{k_2^{2}}\equiv 2f_0\equiv 0\equiv n_0-p_0\hskip
-1mm\pmod{q}.$$

{\bf Case IV:} $q\in \{ 6143, 12289\}$. By $2^6\mid k_i$
$(i=1,2)$, $3\mid k_1$, $3\nmid k_2$ and $u(q)\mid 3\times
2^{12}$, we have $k_1^2\equiv 0\pmod{u(q)}$ and $k_2^2\equiv
2^{12}\pmod{u(q)}$. It follows that
$$f_{k_1^{2}}+f_{k_2^{2}}\equiv f_{2^{12}}\equiv n_0-p_0\pmod{q}.$$

This completes the proof of Theorem \ref{thm2}.
\end{proof}

\begin{proof}[Proof of Theorem \ref{thm1}] Let $M_1=769M$ and
$$n_0\hskip -3mm\pmod{M} \cap 501\hskip -3mm\pmod{769} = n_1\hskip -3mm\pmod{M_1}.$$ Let $n\in
n_1\hskip -1mm\pmod{M_1}$ and $n>0$. Suppose that $n$ can be represented as
 $p+f_{k_1^{2}}
+f_{k_{2}^{2}}$, where $p\in \mathcal{P}$  and $k_1, k_2\in
\mathbb{N}_0$. Since $$n\in n_1\hskip -3mm\pmod{M_1} \subseteq
n_0\hskip -3mm\pmod{M},$$ it follows from Theorem \ref{thm2} that
$p\equiv p_0\hskip -1mm\pmod M$ and $k_i\in \mathcal{K}_i$
$(i=1,2)$ or $k_i\in \mathcal{K}_{3-i}$ $(i=1,2)$. Without loss of
generality, we may assume that $k_i\in \mathcal{K}_i$ $(i=1,2)$.
Then $k_1^2\equiv 0\pmod{3\times 2^{6}}$ and $k_2^2\equiv
64\pmod{3\times 2^{6}}$. Noting that $u(769)=3\times 2^{6}$,
 we have
$$p=n-f_{k_1^{2}}-f_{k_2^{2}}\equiv 501-f_0-f_{64}\equiv
0\hskip -3mm\pmod{769}.$$ It follows that $p=769$. By \eqref{eq1}
and $p\equiv p_0\hskip -1mm\pmod M$, we have $p\equiv 3\pmod 7$.
But $p=769\equiv 6\hskip -1mm\pmod 7$, a contradiction. Hence,
none of integers in $\{ M_1 k+n_1 : k=0,1,\dots \} $ can be
represented as
 $p+f_{k_1^{2}}
+f_{k_{2}^{2}}$, where  $p\in \mathcal{P}$  and $k_1, k_2\in
\mathbb{N}_0$.

This completes the proof of Theorem \ref{thm1}.
\end{proof}

\section{Proof of Theorem \ref{thm3}}

Let $M$, $n_0$, $p_0$, $\mathcal{P}_M$ and $\mathcal{K}_i$
$(i=1,2)$ be as in Section \ref{sec1}.   In this section, we will
employ Theorem \ref{thm2} to prove that there is a positive
proportion of positive integers $n$ with $n\equiv n_0\hskip
-1mm\pmod{M}$ which can be uniquely represented as $p+f_{k_1^{2}}
+f_{k_{2}^{2}}$, where  $p\in \mathcal{P}$  and $k_1, k_2\in
\mathbb{N}_0$ with $k_1\leq k_2$. Let $\alpha=(1+\sqrt{5})/2$ and
$\beta=(1-\sqrt{5})/2$. It is well known that
$$f_n =\frac 1{\sqrt 5} \left( \alpha^n-
\beta^n\right) , \quad n=0,1,\dots .$$

In order to prove Theorem \ref{thm3}, we give the following
lemmas.

\begin{lemma}\label{lem0} If $n\equiv n_0 \hskip -1mm\pmod{M}$, then
$$r(n)=\sharp \{ (p, k_1, k_2 ): n=p+f_{k_1^2}+f_{k_2^2}, p\in
\mathcal{P}_M, k_i\in \mathcal{K}_i (i=1,2)\} .$$
\end{lemma}

\begin{proof} If $r(n)=0$, then the lemma is true trivially. Now let $n\equiv n_0 \hskip -1mm\pmod{M}$ with $r(n)\ge
1$. Then $n$ can be represented as $n=p+f_{k_1^2}+f_{k_2^2}, p\in
\mathcal{P}, k_1, k_2\in \mathbb{N}_0$. By Theorem \ref{thm2},
$p\in \mathcal{P}_M$ and either $k_i\in \mathcal{K}_i (i=1,2)$ or
$k_i\in \mathcal{K}_{3-i} (i=1,2)$. So $k_1\not= k_2$. Therefore,
\begin{align*}r(n)&=\frac 12 \sharp \{ (p, k_1, k_2 ): n=p+f_{k_1^2}+f_{k_2^2}, p\in \mathcal{P}, k_1, k_2\in \mathbb{N}_0\}
\\
&=\frac 12 \sharp \{ (p, k_1, k_2 ): n=p+f_{k_1^2}+f_{k_2^2}, p\in
\mathcal{P}_M, \\
&\hskip 0.8cm k_i\in \mathcal{K}_i (i=1,2) \text{ or }  k_i\in
\mathcal{K}_{3-i} (i=1,2) \}
\\
&=\sharp \{ (p, k_1, k_2 ): n=p+f_{k_1^2}+f_{k_2^2}, p\in
\mathcal{P}_M, k_i\in \mathcal{K}_i (i=1,2)\} .
\end{align*}

This completes the proof of Lemma \ref{lem0}.
\end{proof}

\begin{lemma}\label{lem1} For all sufficiently large numbers $x$, we have
$$\sum_{\substack{n\leq{x}\\ n\equiv n_0\hskip -2mm\pmod{M}}}r(n)=(1+o(1))\frac{x}{2^{11}\times 3^2\varphi(M)\log \alpha},$$
where $\varphi $ is the Euler totient function.
\end{lemma}

\begin{proof}  By Lemma \ref{lem0} and the prime number theorem on arithmetic progressions
(see \cite[Theorem 7.10]{Apostol} or \cite[Chapter
20]{Davenport}), we have
\begin{align*}
\sum_{\substack{n\leq{x}\\ n\equiv n_0 \hskip -2mm\pmod{M}}}r(n)
=& \sharp \{ (p, k_1, k_2 ): p+f_{k_1^2}+f_{k_2^2}\leq x, p\in
\mathcal{P}_M, k_i\in \mathcal{K}_i (i=1,2) \} \\
\leq & \sharp\{p\leq x: p\in \mathcal{P}_M \} \cdot \prod_{i=1}^2 \sharp\{k_i: f_{k_i^2}\leq x, k_i\in \mathcal{K}_i\}\\
=&(1+o(1))\frac{x}{2^{11}\times 3^2 \varphi(M)\log \alpha}
\end{align*}
and
\begin{align*}
\sum_{\substack{n\leq{x}\\ n\equiv n_0 \hskip -2mm\pmod{M}}}r(n)
=&\sharp\{ (p, k_1, k_2 ): p+f_{k_1^2}+f_{k_2^2}\leq x, p\in
\mathcal{P}_M, k_i\in \mathcal{K}_i (i=1,2)\}\\
\geq & \sharp\{p: p\leq x-2{x}/{\log x},  p\in \mathcal{P}_M\}\\
&\cdot \prod_{i=1}^2 \sharp\{k_i: f_{k_i^2}\leq x/\log x, k_i\in \mathcal{K}_i\}\\
=&(1+o(1))\frac{x}{2^{11}\times 3^2 \varphi(M)\log \alpha}.
\end{align*}
This completes the proof of Lemma \ref{lem1}.
\end{proof}

\begin{lemma}\cite[Proposition 2.1]{YDing2022}\label{lem2}
Let $a$ and $K$ be two positive integers, then the number of
solutions of congruent equation $y^2\equiv a \hskip -1mm\pmod{K}$
 does not exceed $4\sqrt{K}$ .
\end{lemma}

The following lemma is a special case of \cite[Theorem 3.12]{Sieve}.

\begin{lemma}\cite[Theorem 3.12]{Sieve}\label{lem3a}
 Let $A$ be a positive number and $k,\ell, m$ positive integers with $(k,\ell
)=1$ and $1\le k\le (\log x)^A$. Then as $x\to +\infty $, we have, uniformly in $k,\ell, m$,  that the number of primes
$p\le x$ with $p\equiv l\hskip -1mm\pmod k$, for which $p+m$ is also a prime,
is less than
$$(8+o_A(1))\prod_{p>2} \left( 1-\frac 1{(p-1)^2}\right) \prod_{\substack{p\mid km\\ p>2}} \frac{p-1}{p-2} \frac{x}{\varphi (k) (\log x)^2}.$$
\end{lemma}

\begin{lemma}\label{lem3}
 Let $A$ be a positive number and $k,\ell, m$ be positive integers with $(k,\ell
)=1$ and $1\le k\le (\log x)^A$. Then as $x\to +\infty $, we have,
uniformly in $k,\ell, m$,  that the number of primes $p\le x$ with
$p\equiv l\hskip -1mm\pmod k$, for which $p+m$ is also a prime, is
less than
$$ (\pi^2+o_A(1))\frac{ x}{\varphi (k) (\log x)^2} \prod_{\substack{p\mid km\\ p>2}} \left(
1+\frac 1{p}\right).$$
\end{lemma}

\begin{proof} Since
\begin{align*}&\prod_{p>2} \left( 1-\frac 1{(p-1)^2}\right) \prod_{\substack{p\mid km\\ p>2}} \frac{(p-1)p}{(p-2)(p+1)}\\
< & \prod_{p>2} \left( 1-\frac 1{(p-1)^2}\right)\frac{(p-1)p}{(p-2)(p+1)}\\
=& \prod_{p>2} \left( 1-\frac 1{p^2}\right)^{-1}
= \sum_{n=0}^\infty \frac 1{(2n+1)^2}=\frac{\pi^2}8,
\end{align*}
Lemma \ref{lem3} follows from Lemma \ref{lem3a}.
\end{proof}

\begin{lemma}\label{lem4} We have
$$\sum\limits_{d=1}^{\infty}\frac{\mu^2(d)}{d\sqrt{[192^2, v(d)]}}<0.23219,$$
where $\mu (d)$ is the M\"obius function.
\end{lemma}

\begin{proof}
The proof begins with the following equality:
\begin{align*}
\sum\limits_{d=1}^{\infty}\frac{\mu^2(d)}{d\sqrt{[192^2,
v(d)]}}&=\sum_{n=1}^\infty \frac
1{\sqrt{[192^2,n]}}\sum_{v(d)=n}\frac{\mu^2(d)}{d}.
\end{align*}
For  $N>192^2$, by Abel's summation we have
\begin{align*}  &\sum_{n=192^2+1}^N \frac 1{\sqrt{n}}\sum_{v(d)=n}\frac{\mu^2(d)}{d}\\
=&\sum_{n=192^2+1}^N \left(\frac 1{\sqrt{n}}-\frac 1{\sqrt{n+1}}\right) \sum_{v(d)\le n}\frac{\mu^2(d)}{d} \\
&+ \frac 1{\sqrt{N+1}}\sum_{v(d)\le N}\frac{\mu^2(d)}{d}-\frac
1{\sqrt{192^2+1}}\sum_{v(d)\le 192^2}\frac{\mu^2(d)}{d}.
\end{align*}
Now we are in a position to estimate
$$ \sum_{v(d)\le n}\frac{\mu^2(d)}{d},\quad n\ge 192^2.$$
Let $d$ be a squarefree integer. Noting that $f_{u(p)}\equiv
f_0\equiv 0\pmod{p}$ for any prime $p$, we know that if $d>1$ and
$v(d)\le n$, then $u(p)\le n$ for every prime divisor $p$ of $d$.
It follows that $d\mid f_1\cdots f_n$. This is also true for
$d=1$. Thus
\begin{align*}
\sum\limits_{v(d)\leq n}\frac{\mu^2(d)}{d}
&\le \sum\limits_{d\mid f_1\cdots f_n}\frac{\mu^2(d)}{d}
=\prod\limits_{p\mid f_1\cdots f_n}\left(1+\frac{1}{p}\right)\\
&<\prod\limits_{p\mid f_1\cdots
f_n}\left(1-\frac{1}{p}\right)^{-1}=\frac{ f_1\cdots f_n}{\varphi(
f_1\cdots f_n)}.
\end{align*}
 By
\cite[(3.41)]{Rosser}, for $k\ge
    3$, we have
    $$\frac{k}{\varphi (k)} \le e^{\gamma} \log \log k+\frac{5}{2 \log \log
        k},$$
where $\gamma$ is the Euler constant. So
\begin{align*}
\sum\limits_{v(d)\leq n}\frac{\mu^2(d)}{d} &<e^{\gamma} \log \log
(f_1\cdots f_n)+\frac{5}{2 \log \log
        (f_1\cdots f_n)}.
\end{align*}
Noting that  $f_k<\alpha^k$ for all $k\ge 1$, we have
$f_1\cdots f_n<\alpha^{1+2+\cdots +n}\le \alpha ^{n^2}$. It follows that
$$\log \log (f_1\cdots f_n)<\log\log \alpha ^{n^2}=2\log n+\log\log \alpha .$$
Since $e^{\gamma}\log\log \alpha <-1.3$ and
$$\frac{5}{2 \log \log
        (f_1\cdots f_n)}<\frac{5}{2 \log \log
        (f_1\cdots f_{7})}<1.2,\quad n\ge 7,$$
it follows that
$$\sum\limits_{v(d)\leq n}\frac{\mu^2(d)}{d}<2e^{\gamma}\log n,\quad n\ge 7.$$
Hence
\begin{align*}
&\sum_{n=1}^\infty \frac 1{\sqrt{[192^2,n]}}\sum_{v(d)=n}\frac{\mu^2(d)}{d}\\
=&\sum_{n=1}^{192^2} \frac 1{\sqrt{[192^2,n]}}\sum_{v(d)=n}\frac{\mu^2(d)}{d}+\sum_{n=192^2+1}^\infty \frac 1{\sqrt{[192^2,n]}}\sum_{v(d)=n}\frac{\mu^2(d)}{d}\\
\le & \sum_{n=1}^{192^2} \frac 1{192}\sum_{v(d)=n}\frac{\mu^2(d)}{d}+\sum_{n=192^2+1}^\infty \frac 1{\sqrt{n}}\sum_{v(d)=n}\frac{\mu^2(d)}{d}\\
\le &  \frac 1{192}\sum_{v(d)\le 192^2}\frac{\mu^2(d)}{d}+\sum_{n=192^2+1}^\infty \left(\frac 1{\sqrt{n}}-\frac 1{\sqrt{n+1}}\right)\sum_{v(d)\le n}\frac{\mu^2(d)}{d}\\
& -\frac 1{\sqrt{192^2+1}}\sum_{v(d)\le 192^2}\frac{\mu^2(d)}{d}\\
=&\left(\frac 1{192}-\frac 1{\sqrt{192^2+1}}\right) \sum_{v(d)\le 192^2}\frac{\mu^2(d)}{d} +\sum_{n=192^2+1}^\infty \left(\frac 1{\sqrt{n}}-\frac 1{\sqrt{n+1}}\right)\sum_{v(d)\le n}\frac{\mu^2(d)}{d} \\
\le &2e^{\gamma} \left(\frac 1{192}-\frac 1{\sqrt{192^2+1}}\right) \log 192^2+2e^{\gamma}\sum_{n=192^2+1}^\infty \left(\frac 1{\sqrt{n}}-\frac 1{\sqrt{n+1}}\right)\log n\\
= &2e^{\gamma} \left(\frac 1{192}-\frac 1{\sqrt{192^2+1}}\right)
\log 192^2
+e^{\gamma}\sum_{n=192^2+1}^\infty \int_{n}^{n+1} x^{-\frac 32} \log n dx \\
\le &2e^{\gamma} \left(\frac 1{192}-\frac 1{\sqrt{192^2+1}}\right)
\log 192^2+e^{\gamma} \int_{192^2+1}^\infty x^{-\frac 32}\log x dx
<0.23219.
\end{align*}

This completes the proof of Lemma \ref{lem4}.
\end{proof}

\begin{lemma}\label{lem5} For all sufficiently large numbers $x$, we have
$$\sum_{\substack{n\leq{x}\\ n\equiv n_0\hskip -2mm\pmod{M}}}r(n)(r(n)-1)
\leq \frac{0.2322 \pi^2}{2^{11}\times 3^4\varphi (M) (\log \alpha
)^2} x.$$
\end{lemma}

\begin{proof} Let $x$ be a sufficiently large number. If $f_{k_1^{2}}+f_{k_{2}^{2}} =f_{l_1^{2}}+f_{l_{2}^{2}}$
with $k_i, l_i\in \mathcal{K}_i$  $(i=1,2)$, then by
$f_m=f_{m-1}+f_{m-2}$ for all $m\ge 2$, we have $k_i=l_i$
$(i=1,2)$. In view of Lemma \ref{lem0} and Lemma \ref{lem3}, we
have
\begin{align*}
&\sum_{\substack{n\leq{x}\\ n\equiv n_0\hskip -2mm\pmod{M}}}r(n)(r(n)-1)\\
=&\# \{ (p_1, p_2, k_1, k_2, l_1, l_2) :
p_1+f_{k_1^{2}}+f_{k_{2}^{2}}
=p_2+f_{l_1^{2}}+f_{l_{2}^{2}}\le x,\\
&~~ p_i\equiv p_0\hskip -3mm\pmod{M}, p_i\in \mathcal{P} (i=1,2), p_1\not= p_2, k_i, l_i\in \mathcal{K}_i (i=1,2)\}\\
= &\sum_{\substack{k_i, l_i\in \mathcal{K}_i\\
f_{k_i^{2}}\le x, f_{l_i^{2}}\le x (i=1,2)\\
f_{k_1^{2}}+f_{k_{2}^{2}} \not= f_{l_1^{2}}+f_{l_{2}^{2}}}}
\sum_{\substack{ p_1+f_{k_1^{2}}+f_{k_{2}^{2}}
=p_2+f_{l_1^{2}}+f_{l_{2}^{2}}\le x\\ p_i\equiv p_0\hskip -2mm\pmod{M}, p_i\in \mathcal{P} (i=1,2)}} 1\\
\le & \sum_{\substack{k_i, l_i\in \mathcal{K}_i\\ f_{k_i^{2}}\le
x, f_{l_i^{2}}\le x (i=1,2)\\ f_{k_1^{2}}+f_{k_{2}^{2}} \not=
f_{l_1^{2}}+f_{l_{2}^{2}}}}
\sum_{\substack{ p_1\le x\\ p_1\equiv p_0\hskip -2mm\pmod{M}, p_1\in \mathcal{P} \\ p_1+f_{k_1^{2}}+f_{k_{2}^{2}}-f_{l_1^{2}}-f_{l_{2}^{2}}\in \mathcal{P}}} 1\\
\le & \sum_{\substack{k_i, l_i\in \mathcal{K}_i\\ f_{k_i^{2}}\le
x, f_{l_i^{2}}\le x (i=1,2)\\ f_{k_1^{2}}+f_{k_{2}^{2}} \not=
f_{l_1^{2}}+f_{l_{2}^{2}}}} \frac{(\pi^2 +o(1)) x}{\varphi (M)
(\log x)^2} \prod_{\substack{p\mid
(f_{k_1^{2}}+f_{k_{2}^{2}}-f_{l_1^{2}}-f_{l_{2}^{2}})M\\ p>2}}
\left(
1+\frac 1{p}\right) \\
\le & \frac{2(\pi^2 +o(1)) x}{ 3\varphi (M) (\log x)^2}
\prod_{\substack{p\mid M}}\left( 1+\frac 1{p}\right)
\sum_{\substack{k_i, l_i\in \mathcal{K}_i\\ f_{k_i^{2}}\le x,
f_{l_i^{2}}\le x (i=1,2)\\ f_{k_1^{2}}+f_{k_{2}^{2}} \not=
f_{l_1^{2}}+f_{l_{2}^{2}}}} \prod_{\substack{p\mid
f_{k_1^{2}}+f_{k_{2}^{2}}-f_{l_1^{2}}-f_{l_{2}^{2}}\\ p\nmid M}}
\left( 1+\frac 1{p}\right)  .\end{align*} Given $k_i, l_i\in
\mathcal{K}_i$ $(i=1,2)$ such that $$f_{k_i^{2}}\le x,
f_{l_i^{2}}\le x\ (i=1,2), f_{k_1^{2}}+f_{k_{2}^{2}} \not=
f_{l_1^{2}}+f_{l_{2}^{2}}.$$ Then $1\le
|f_{k_1^{2}}+f_{k_{2}^{2}}-f_{l_1^{2}}-f_{l_{2}^{2}}|\le 2x$. The
number of prime divisors $p$ of $
f_{k_1^{2}}+f_{k_{2}^{2}}-f_{l_1^{2}}-f_{l_{2}^{2}}$ with $p> \log
x$ is less than
$$\frac{\log 2x}{\log\log x}.$$
It follows that
$$\prod_{\substack{p\mid f_{k_1^{2}}+f_{k_{2}^{2}}-f_{l_1^{2}}-f_{l_{2}^{2}}\\ p>\log x}} \left(
1+\frac 1{p}\right)  <\left(
1+\frac 1{\log x}\right)^{\log (2x)/\log\log x}=1+o(1).$$
Hence
\begin{align*}&\sum_{\substack{n\leq{x}\\ n\equiv n_0\hskip -2mm\pmod{M}}}r(n)(r(n)-1)\\
\le & \frac{2(\pi^2 +o(1)) x}{ 3\varphi (M) (\log
x)^2}\prod_{p\mid M}\left( 1+\frac 1{p}\right)
\sum_{\substack{k_i, l_i\in \mathcal{K}_i\\ f_{k_i^{2}}\le x,
f_{l_i^{2}}\le x (i=1,2)\\ f_{k_1^{2}}+f_{k_{2}^{2}} \not=
f_{l_1^{2}}+f_{l_{2}^{2}}}} \prod_{\substack{p\mid
f_{k_1^{2}}+f_{k_{2}^{2}}-f_{l_1^{2}}-f_{l_{2}^{2}}\\  p\le \log
x, p\nmid M}} \left(
1+\frac 1{p}\right) \\
= & \frac{2(\pi^2 +o(1)) x}{ 3\varphi (M) (\log x)^2}\prod_{p\mid
M}\left( 1+\frac 1{p}\right)  \sum_{\substack{k_i, l_i\in
\mathcal{K}_i\\ f_{k_i^{2}}\le x, f_{l_i^{2}}\le x (i=1,2)\\
f_{k_1^{2}}+f_{k_{2}^{2}} \not= f_{l_1^{2}}+f_{l_{2}^{2}}}}
\sum_{\substack{d\mid
f_{k_1^{2}}+f_{k_{2}^{2}}-f_{l_1^{2}}-f_{l_{2}^{2}}\\ (d, M)=1,
P(d)\le \log x}}
\frac{\mu (d)^2}{d}\\
= & \frac{2(\pi^2 +o(1)) x}{ 3\varphi (M) (\log x)^2} \prod_{p\mid
M} \left(
1+\frac 1{p}\right) \sum_{\substack{k_2\in \mathcal{K}_2, f_{k_2^{2}}\le x\\
l_i\in \mathcal{K}_i, f_{l_i^{2}}\le x (i=1,2)}}\sum_{\substack{k_1\in \mathcal{K}_1, f_{k_1^{2}}\le x\\ f_{k_1^{2}}+f_{k_{2}^{2}}
\not= f_{l_1^{2}}+f_{l_{2}^{2}}}} \sum_{\substack{d\mid f_{k_1^{2}}+f_{k_{2}^{2}}-f_{l_1^{2}}-f_{l_{2}^{2}}\\ (d,M)=1, P(d)\le \log x}} \frac{\mu (d)^2}{d}\\
= & \frac{2(\pi^2 +o(1)) x}{ 3\varphi (M) (\log x)^2} \prod_{p\mid
M} \left(
1+\frac 1{p}\right) \sum_{\substack{k_2\in \mathcal{K}_2, f_{k_2^{2}}\le x\\
l_i\in \mathcal{K}_i, f_{l_i^{2}}\le x (i=1,2)}}\sum_{\substack{d=1\\ (d,M)=1\\ P(d)\le \log x}}^\infty \frac{\mu (d)^2}{d} \sum_{\substack{k_1\in \mathcal{K}_1, f_{k_1^{2}}\le x\\ f_{k_1^{2}}+f_{k_{2}^{2}}
\not= f_{l_1^{2}}+f_{l_{2}^{2}}\\ d\mid f_{k_1^{2}}+f_{k_{2}^{2}}-f_{l_1^{2}}-f_{l_{2}^{2}}}}1\\
\le & \frac{2(\pi^2 +o(1)) x}{ 3\varphi (M) (\log
x)^2}\prod_{p\mid M} \left(
1+\frac 1{p}\right) \sum_{\substack{k_2\in \mathcal{K}_2, f_{k_2^{2}}\le x\\
l_i\in \mathcal{K}_i, f_{l_i^{2}}\le x (i=1,2)}}\sum_{\substack{d=1\\ (d,M)=1\\ P(d)\le \log x}}^\infty
\frac{\mu (d)^2}{d} \sum_{\substack{k_1\in \mathcal{K}_1, k_1\le L\\ d\mid f_{k_1^{2}}+f_{k_{2}^{2}}-f_{l_1^{2}}-f_{l_{2}^{2}}}}1,
\end{align*}
where $P(1)=0$,  $P(d)$ denotes the largest prime divisor of $d$
for $d>1$ and $L=\sqrt{\log (\sqrt 5 x +1)/\log \alpha}$. Now we
are in a position to estimate the upper bound of
\begin{equation}\label{e1}
 \sum_{\substack{k_1\in \mathcal{K}_1, k_1\le L\\ d\mid f_{k_1^{2}}+f_{k_{2}^{2}}-f_{l_1^{2}}-f_{l_{2}^{2}}}}1. \end{equation}
Given $k_2, l_1, l_2, d$. If there is at most one $k_1\leq L$ with
$k_1\in \mathcal{K}_1$ such that $d\mid
f_{k_1^{2}}+f_{k_{2}^{2}}-f_{l_1^{2}}-f_{l_{2}^{2}}$, then
\eqref{e1} does not exceed 1. Now we assume that there are at
least two such $k_1$. Let $h_1$ be the least such $k_1$. Since
$d\mid f_{k_1^{2}}+f_{k_{2}^{2}}-f_{l_1^{2}}-f_{l_{2}^{2}}$ and
$d\mid f_{h_1^{2}}+f_{k_{2}^{2}}-f_{l_1^{2}}-f_{l_{2}^{2}}$, it
follows that $d \mid f_{k_{1}^2}-f_{h_{1}^2}$. We assume that
$d>1$. By the definition of $v(d)$, there exists a prime divisor
$p_d$ of $d$ such that $v(d)=u(p_d)$. So $p_d \mid
f_{k_{1}^2}-f_{h_{1}^2}$. Let $B_d$ be the set of all integers
$0\le l<v(d)$ such that $f_l\equiv f_{h_{1}^2}\hskip
-1mm\pmod{p_d}$. In view of \cite{Schinzel} or \cite{Somer},
$|B_d|\le 4$. Thus, $p_d \mid f_{k_{1}^2}-f_{l}$, $l\in B_d$. It
follows that
$$k_{1}^2\equiv l\hskip -3mm\pmod{v(d)},\quad k_1\in \mathcal{K}_1, l\in B_d.$$
Noting that $k_1\in \mathcal{K}_1$, we have $192\mid k_1$. So
$k_1^2\equiv 0\pmod{192^2}$. Let $C_d$ be set of all integers
$0\le c<[192^2, v(d)]$ such that for some $l\in B_d$,
$$c\equiv l\hskip -3mm\pmod{v(d)}, \quad c\equiv 0\hskip -3mm\pmod{192^2}.$$
Then $|C_d|=|B_d|\le 4$ and
\begin{equation}\label{e5}k_{1}^2\equiv c\hskip -3mm\pmod{[192^2, v(d)]},\quad k_{1}\le L, c\in C_d.\end{equation}
Hence, \eqref{e1} does not exceed the number of solutions of
\eqref{e5}. If $L^{\frac{4}{3}}<[192^2, v(d)]$, then $0\le
k_{1}^2- c \le L^2<L^{\frac{2}{3}}[192^2, v(d)]$. Thus, the number
of positive integers $k_1$ with \eqref{e5} does not exceed
$4(L^{\frac{2}{3}}+1)$. If $L^{\frac{4}{3}}>[192^2, v(d)]$, then
by Lemma \ref{lem2}, the number of positive integers $k_1$ with
\eqref{e5} does not exceed
\begin{align*}
16\sqrt{[192^2, v(d)]}\left(\frac{L}{[192^2, v(d)]}+1\right)
&=\frac{16L}{\sqrt{[192^2, v(d)]}}+16\sqrt{[192^2, v(d)]}\\
&\leq \frac{16L}{\sqrt{[192^2, v(d)]}} +16L^{\frac{2}{3}}.
\end{align*}
In all cases, if $d>1$, then the number of solutions of \eqref{e5}
does not exceed $$\frac{16L}{\sqrt{[192^2, v(d)]}}
+16L^{\frac{2}{3}}$$ and so does \eqref{e1}. If $d=1$, then by
$v(1)=1$, \eqref{e1} is
$$\sum_{k_1\in \mathcal{K}_1, k_1\le L}1\le \frac{L}{192}+1\le \frac{16L}{\sqrt{[192^2, v(d)]}} +16L^{\frac{2}{3}}.
$$
Hence
\begin{align*}
\sum_{\substack{d=1\\ (d,M)=1\\ P(d)\le \log x}}^\infty \frac{\mu
(d)^2}{d} \sum_{\substack{k_1\in \mathcal{K}_1, k_1\le L\\ d\mid
f_{k_1^{2}}+f_{k_{2}^{2}}-f_{l_1^{2}}-f_{l_{2}^{2}}}}1 \le &
\sum_{\substack{d=1\\ (d,M)=1\\ P(d)\le \log x}}^\infty \frac{\mu
(d)^2}{d} \left( \frac{16L}{\sqrt{[192^2, v(d)]}}
+16L^{\frac{2}{3}}\right).
\end{align*}
By Mertens' formula, we have \begin{equation}\label{eqa1}\sum_{\substack{d=1\\ (d,M)=1\\
P(d)\le \log x}}^\infty \frac{\mu (d)^2}{d} \le \sum_{\substack{
d=1\\  P(d)\le \log x}}^\infty \frac{\mu (d)^2}{d} =\prod_{p\le
\log x} \left( 1+\frac 1p\right) =O(\log\log x).\end{equation} It
follows that
\begin{align*}
\sum_{\substack{d=1\\ (d,M)=1\\ P(d)\le \log x}}^\infty \frac{\mu
(d)^2}{d} \sum_{\substack{k_1\in \mathcal{K}_1, k_1\le L\\ d\mid
f_{k_1^{2}}+f_{k_{2}^{2}}-f_{l_1^{2}}-f_{l_{2}^{2}}}}1 \le
&\sum_{\substack{d=1\\ (d,M)=1}}^\infty \frac{16L\mu
(d)^2}{d\sqrt{[192^2, v(d)]}} +o(L).
\end{align*}
Thus
\begin{align*}
&\sum_{\substack{n\leq{x}\\ n\equiv n_0\hskip -2mm\pmod{M}}}r(n)(r(n)-1)\\
\le & \frac{2(\pi^2 +o(1)) x}{ 3\varphi (M) (\log x)^2}
\prod_{p\mid M} \left(
1+\frac 1{p}\right) \cdot \frac{2^2L^3}{192^3}\left( \sum_{\substack{d=1\\ (d,M)=1}}^\infty \frac{16L\mu (d)^2}{d\sqrt{[192^2, v(d)]}} +o(L)\right)\\
\le & \frac{(\pi^2 +o(1)) x}{2^{11}\times 3^4 \varphi (M) (\log
\alpha )^2}
\prod_{p\mid M} \left( 1+\frac 1{p}\right) \sum_{\substack{d=1\\
(d,M)=1}}^\infty \frac{\mu (d)^2}{d\sqrt{[192^2, v(d)]}} +o(x).
\end{align*}

Let $d_1, d_2$ be positive integers with
$d_1\mid M$ and $(d_2, M)=1$. By \eqref{e4}, $v(d_1)\mid
192^2$. Noting that $$v(d_1d_2)=\max \{ v(d_1), v(d_2)\} ,$$ we
have $[192^2, v(d_1d_2)]\le [192^2, v(d_2)]$. Thus, by Lemma
\ref{lem4}, we have
\begin{align*}
&\prod_{\substack{p\mid M}} \left(
1+\frac 1{p}\right)\sum_{\substack{d=1\\
(d, M)=1}}^\infty\frac{\mu^2(d)}{d\sqrt{[192^2, v(d)]}}\\
&\leq \sum_{\substack{d_1\mid
M}}\frac{\mu^2(d_1)}{d_1}
 \sum _{\substack{d_2=1\\(d_2, M)=1}}^{\infty}\frac{\mu^2(d_2)}{d_2\sqrt{[192^2, v(d_2)]}}\\
&=   \sum _{\substack{d_2=1\\(d_2, M)=1\\ d_1\mid M}}^{\infty}\frac{\mu^2(d_1d_2)}{d_1d_2\sqrt{[192^2, v(d_2)]}}\\
&\leq  \sum _{\substack{d_2=1\\(d_2, M)=1\\ d_1\mid M}}^{\infty}\frac{\mu^2(d_1d_2)}{d_1d_2\sqrt{[192^2, v(d_1d_2)]}}\\
&= \sum
_{\substack{d=1}}^{\infty}\frac{\mu^2(d)}{d\sqrt{[192^2,
v(d)]}}\\
&<0.23219.
\end{align*}
It follows that
\begin{align*}
&\sum_{\substack{n\leq{x}\\ n\equiv n_0\hskip -2mm
\pmod{M}}}r(n)(r(n)-1) \le \frac{0.2322 \pi^2}{2^{11}\times
3^4\varphi (M) (\log \alpha )^2} x.
\end{align*}
This completes the proof of Lemma \ref{lem5}.
\end{proof}

\begin{proof}[Proof of Theorem \ref{thm3}]
By Lemma \ref{lem1} and Lemma \ref{lem5}, for all sufficiently
large $x$, we have
\begin{align*}
\sum_{\substack{n\leq{x}\\ n\equiv n_0\hskip
-2mm\pmod{M}\\r(n)=1}}1 =&\sum_{\substack{n\leq{x}\\ n\equiv
n_0\hskip -2mm\pmod{M}}}r(n)-
\sum_{\substack{n\leq{x}\\ n\equiv n_0\hskip -2mm\pmod{M}\\r(n)\geq 2}}r(n)\\
\ge &\sum_{\substack{n\leq{x}\\ n\equiv n_0\hskip -2mm\pmod{M}}}r(n)-
\sum_{\substack{n\leq{x}\\ n\equiv n_0\hskip -2mm\pmod{M}\\r(n)\geq 2}}r(n)(r(n)-1)\\
\geq & (1+o(1))\frac{x}{2^{11}\times 3^2\varphi(M)\log \alpha}-\frac{0.2322 \pi^2}{2^{11}\times 3^4\varphi (M)(\log \alpha)^2}\, x\\
= & \frac{9\log \alpha-0.2322 \pi^2 +o(1)}{2^{11}\times 3^4\varphi(M)(\log \alpha)^2}\, x\\
> &\frac{2}{2^{11}\times 3^4\varphi(M)(\log \alpha)^2}\, x.
\end{align*}
This completes the proof of Theorem \ref{thm3}.

\end{proof}

\section{Proof of Theorem \ref{thm4}}

Recall that $\alpha =(1+\sqrt 5)/2$ and
$$r(n)=\sharp \{ (p, k_1, k_2 ): n=p+f_{k_1^2}+f_{k_2^2}, p\in \mathcal{P}, k_1, k_2\in \mathbb{N}_0, k_1\le
k_2\} .$$

\begin{lemma}\label{lema1}We have
$$\sum_{n\le x} r(n)= (1+o(1))\frac{x}{2\log \alpha}.$$
\end{lemma}

\begin{proof} In the following, for simplicity we always assume that $p\in \mathcal{P}$ and $k_1, k_2\in
\mathbb{N}_0$. By the prime number theorem, we have
\begin{align*}
\sum_{n\le x} r(n) =& \sharp \{ (p, k_1, k_2 ):
p+f_{k_1^2}+f_{k_2^2}\leq x, k_1\le
k_2 \} \\
\leq & \sharp\{p\leq x \} \cdot  \sharp\{ (k_1, k_2): f_{k_1^2}\le f_{k_2^2} \leq x\}\\
=&(1+o(1))\frac{x}{2\log \alpha}
\end{align*}
and
\begin{align*}
\sum_{n\le x} r(n) =&\sharp\{ (p, k_1, k_2 ):
p+f_{k_1^2}+f_{k_2^2}\leq x,  k_1\le
k_2 \}\\
\geq & \sharp\{ p\leq x-2{x}/{\log x}\}
\cdot \sharp\{ (k_1, k_2): f_{k_1^2}\le f_{k_2^2} \leq {x}/{\log x}\}\\
=&(1+o(1))\frac{x}{2\log \alpha}.
\end{align*}

This completes the proof of Lemma \ref{lema1}.
\end{proof}

\begin{lemma}\label{lema2}We have
$$\sum_{n\le x} r(n)^2\ll x.$$
\end{lemma}

\begin{proof} In the following, for simplicity we always assume that $p_1,p_2\in \mathcal{P}$, $k_1, k_2, l_1, l_2\in
\mathbb{N}_0$, $k_1\le k_2$ and $l_1\le l_2$. We follow the proof
of Lemma \ref{lem5}. By the definition of $r(n)$ and Lemma
\ref{lem3},
\begin{align}\label{u0}
\sum_{n\le x} r(n)^2 &= \sharp \{ (p_1, k_1, k_2, p_2, l_1, l_2) :
p_1+f_{k_1^2}+f_{k_2^2}= p_2+f_{l_1^2}+f_{l_2^2}\leq
x\} \nonumber\\
&\le  \sharp \{ (p_1, k_1, k_2, p_2, l_1, l_2) :
p_1+f_{k_1^2}+f_{k_2^2}= p_1+f_{l_1^2}+f_{l_2^2}\leq
x, p_1=p_2\} \nonumber\\
&\quad +\sharp \{ (p_1, k_1, k_2, p_2, l_1, l_2) :
p_1+f_{k_1^2}+f_{k_2^2}= p_2+f_{l_1^2}+f_{l_2^2}\leq
x,p_1\not= p_2 \} \nonumber\\
&\ll x+ \sum_{\substack{f_{k_1^2}\le x,f_{k_2^2}\le x\\
f_{l_1^2}+f_{l_2^2}\leq x\\ f_{k_1^2}+f_{k_2^2}\not =
f_{l_1^2}+f_{l_2^2}}}
 \sharp \{ p_1 : p_1\le x, p_1+f_{k_1^2}+f_{k_2^2}-f_{l_1^2}-f_{l_2^2}\in \mathcal{P}\} \nonumber\\
&\ll x+ \sum_{\substack{f_{k_2^2}\le x\\ f_{l_1^2}+f_{l_2^2}\leq
x}}
 \sum_{\substack{f_{k_1^2}\le x\\ f_{k_1^2}+f_{k_2^2}\not = f_{l_1^2}+f_{l_2^2}}}
\frac{x}{(\log x)^2} \prod_{ p\mid
f_{k_1^2}+f_{k_2^2}-f_{l_1^2}-f_{l_2^2}} \left( 1+\frac
1{p}\right).
\end{align}
Given $k_2, l_1, l_2$ with $f_{k_2^2}\le x$ and $
f_{l_1^2}+f_{l_2^2}\leq x$. As the arguments of Lemma \ref{lem5},
we have
\begin{align}\label{u1}
&\sum_{\substack{f_{k_1^2}\le x\\ f_{k_1^2}+f_{k_2^2}\not =
f_{l_1^2}+f_{l_2^2}}} \prod_{ p\mid
f_{k_1^2}+f_{k_2^2}-f_{l_1^2}-f_{l_2^2}} \left( 1+\frac
1{p}\right)\nonumber\\
& = (1+o(1)) \sum_{\substack{f_{k_1^2}\le x\\
f_{k_1^2}+f_{k_2^2}\not = f_{l_1^2}+f_{l_2^2}}}
\prod_{\substack{p\mid f_{k_1^2}+f_{k_2^2}-f_{l_1^2}-f_{l_2^2}\\
p\le \log x}} \left( 1+\frac
1{p}\right)\nonumber\\
& = (1+o(1)) \sum_{\substack{f_{k_1^2}\le x\\
f_{k_1^2}+f_{k_2^2}\not = f_{l_1^2}+f_{l_2^2}}}\sum_{
\substack{d\mid
f_{k_1^2}+f_{k_2^2}-f_{l_1^2}-f_{l_2^2}\\ P(d)\le \log x}} \frac{\mu(d)^2}{d}\nonumber\\
&\ll \sum_{d=1, P(d)\le \log x}^\infty \frac{\mu(d)^2}{d} \sum_{
\substack{f_{k_1^2}\le x\\ d\mid
f_{k_1^2}+f_{k_2^2}-f_{l_1^2}-f_{l_2^2}}}  1\nonumber\\
&\ll \sum_{d=1, P(d)\le \log x}^\infty \frac{\mu(d)^2}{d} \left(
\frac{L}{\sqrt{v(d)}}+L^{\frac 23}\right).
\end{align}
In view of Lemma \ref{lem4},
\begin{equation}\label{eqa2}\sum_{d=1}^\infty\frac{\mu(d)^2}{d\sqrt{v(d)}}\le
192\sum\limits_{d=1}^{\infty}\frac{\mu^2(d)}{d\sqrt{[192^2,
v(d)]}}<192\times 0.23219.\end{equation} By \eqref{eqa1},
\eqref{u1} and \eqref{eqa2},
$$\sum_{\substack{f_{k_1^2}\le x\\ f_{k_1^2}+f_{k_2^2}\not =
f_{l_1^2}+f_{l_2^2}}} \prod_{ p\mid
f_{k_1^2}+f_{k_2^2}-f_{l_1^2}-f_{l_2^2}} \left( 1+\frac
1{p}\right)\ll L.$$ It follows from \eqref{u0} that
$$\sum_{n\le x} r(n)^2\ll x+ \sum_{\substack{f_{k_2^2}\le x\\ f_{l_1^2}+f_{l_2^2}\leq
x}} \frac{x}{(\log x)^2}L \ll x+ \frac{x}{(\log x)^2} L^4\ll x.$$

This completes the proof of Lemma \ref{lema2}.
\end{proof}

\begin{proof}[Proof of Theorem \ref{thm4}] Noting that $2\log
\alpha <1$, by Lemma \ref{lema1}, we have
$$\sum_{\substack{n\le x\\ r(n)\ge 2}}r(n)=(1+o(1))\frac{x}{2\log
\alpha} -\sum_{\substack{n\le x\\ r(n)=1}}1 \gg x.$$ By the
Cauchy-Schwarz inequality and Lemma \ref{lema2}, we have
\begin{align*}x^2&\ll \Big( \sum_{\substack{n\le x\\ r(n)\ge 2}}
r(n) \Big)^2\le \Big( \sum_{\substack{n\le x\\ r(n)\ge 2}} 1 \Big)
\Big( \sum_{\substack{n\le x\\ r(n)\ge 2}} r(n)^2
\Big)\\
&\ll \Big( \sum_{\substack{n\le x\\ r(n)\ge 2}} 1 \Big) \Big(
\sum_{n\le x} r(n)^2 \Big) \ll x \Big( \sum_{\substack{n\le x\\
r(n)\ge 2}} 1 \Big) .
\end{align*}
It follows that
$$\sum_{\substack{n\le x\\ r(n)\ge 2}} 1\gg x.$$

This completes the proof of Theorem \ref{thm4}.
\end{proof}

\section*{Acknowledgments}

This work was supported by the National Natural Science Foundation
of China, Grant No. 12171243.

\end{document}